\documentclass[12pt, a4paper]{article}

\usepackage{amsmath,amssymb}
\usepackage{times,fancyhdr}
\usepackage{graphicx}
\usepackage{amsfonts}

\newtheorem{theorem}{Theorem}
\newtheorem{algorithm}{Algorithm}

\begin{document}
\begin{center}
{\bf \Large The LBFGS Quasi-Newtonian Method for Molecular Modeling Prion AGAAAAGA Amyloid Fibrils}\\

\vskip 0.3cm

{\bf {\small Jiapu Zhang$^1$, Yating Hou$^2$, Yiju Wang$^{2,3}$, Changyu Wang$^{2,4,5}$, Xiangsun Zhang$^{4}$}}\\

\vskip 0.3cm

{\small
$^1$Centre for Informatics and Applied Optimization \&\\
Graduate School of Sciences, Information Technology and Engineering,\\
The University of Ballarat, MT Helen Campus, Victoria 3353, Australia,\\
E-mail: j.zhang@ballarat.edu.au. Tel.: +61 4 2348 7360, +61 3 5327 9084.\\ 
$^2$School of Management Science, Qufu Normal University,\\
Rizhao Campus, Rizhao, Shandong 276800, P. R. China, \&\\
$^3$Department of Applied Mathematics, Hong Kong Polytechnic University,\\
Hung Hum, Kowloon, Hong Kong, P. R. China.\\
$^4$Institute of Applied Mathematics, Academia Sinica, Beijing 100080, P. R. China.\\
$^5$Dept. of Applied Maths, Dalian University of Technology, Dalian 116024, P. R. China.}
\end{center}

\begin{abstract}
Experimental X-ray crystallography, NMR (Nuclear Magnetic Resonance) spectroscopy, dual polarization interferometry, etc are indeed very powerful tools to determine the 3-Dimensional structure of a protein (including the membrane protein); theoretical mathematical and physical computational approaches can also allow us to obtain a description of the protein 3D structure at a submicroscopic level for some unstable, noncrystalline and insoluble proteins. X-ray crystallography finds the X-ray final structure of a protein, which usually need refinements using theoretical protocols in order to produce a better structure. This means theoretical methods are also important in determinations of protein structures. Optimization is always needed in the computer-aided drug design, structure-based drug design, molecular dynamics, and quantum and molecular mechanics. This paper introduces some optimization algorithms used in these research fields  and presents a new theoretical computational method - an improved LBFGS Quasi-Newtonian mathematical optimization method - to produce 3D structures of Prion AGAAAAGA amyloid fibrils (which are unstable, noncrystalline and insoluble), from the potential energy minimization point of view. Because the NMR or X-ray structure of the hydrophobic region AGAAAAGA of prion proteins has not yet been determined, the model constructed by this paper can be used as a reference for experimental studies on this region, and may be useful in furthering the goals of medicinal chemistry in this field.
\end{abstract}


\noindent {\bf Highlights:} 
$\blacktriangleright$ Determinations of 3D structures of unstable, noncrystalline and insoluble proteins.
$\blacktriangleright$ Theoretical computer, mathematical/physical computational approaches and concepts.
$\blacktriangleright$ Time-consuming, costly X-ray crystallography and NMR spectroscopy powerful tools.
$\blacktriangleright$ Mathematical local search optimization methods to solve three-body physical problems. 
$\blacktriangleright$ Using computational chemistry Amber/Gromacs packages to refine molecular models.

\section{Introduction}
Neurodegenerative diseases including Parkinson's, Alzheimer's, Huntington's, and Prion's were found they all featured amyloid fibrils \cite{chiang_etal2008, irvine_etal2008, ferreira_etal2007, nature_editorial2001, truant_etal2008, weydt_etal2006}. Amyloid is characterized by a cross-$\beta$ sheet quaternary structure and recent X-ray diffraction studies of microcrystals revealed atomistic details of core region of amyloid \cite{nelson_etal2005, sawaya_etal2007}. All the quaternary structures of amyloid cross-$\beta$ spines can be reduced to one of the 8 classes of steric zippers of \cite{sawaya_etal2007}, with strong van der Waals (vdW) interactions between $\beta$-sheets and hydrogen bonds (HBs) to maintain the $\beta$-strands. A new era in the structural analysis of amyloids started from the `steric zipper'- $\beta$-sheets \cite{nelson_etal2005}. As the two $\beta$-sheets zip up, Hydrophobic Packings (HPs \& vdWs) have been formed. The extension of the `steric zipper' above and below (i.e. the $\beta$-strands) is maintained by Hydrogen Bonds (HBs) (but usually there is no HB between the two $\beta$-sheets). This is the common structure associated with some 20 neurodegenerative amyloid diseases.\\

We first do some mathematical analysis for the common structure. Let $r$ be the distance between two atoms, the vdW contacts of the two atoms are described by the Lennard-Jones (LJ) potential energy:
$$ V_{LJ} (r) = 4 \epsilon \left( (\frac{\sigma }{r} )^{12} -(\frac{\sigma }{r} )^6 \right),\eqno(1.1)$$
where $\epsilon$ is the depth of the potential well and $\sigma$ is the atom diameter; these parameters can be fitted to reproduce experimental data or deduced from results of accurate quantum chemistry calculations. The $(\frac{\sigma }{r} )^{12}$ item describes repulsion and the $-(\frac{\sigma }{r} )^6$ item describes attraction (Figure 1).  
\begin{figure}[h!] \label{Fig01-LJ_potential}
\includegraphics[width=5.0in]{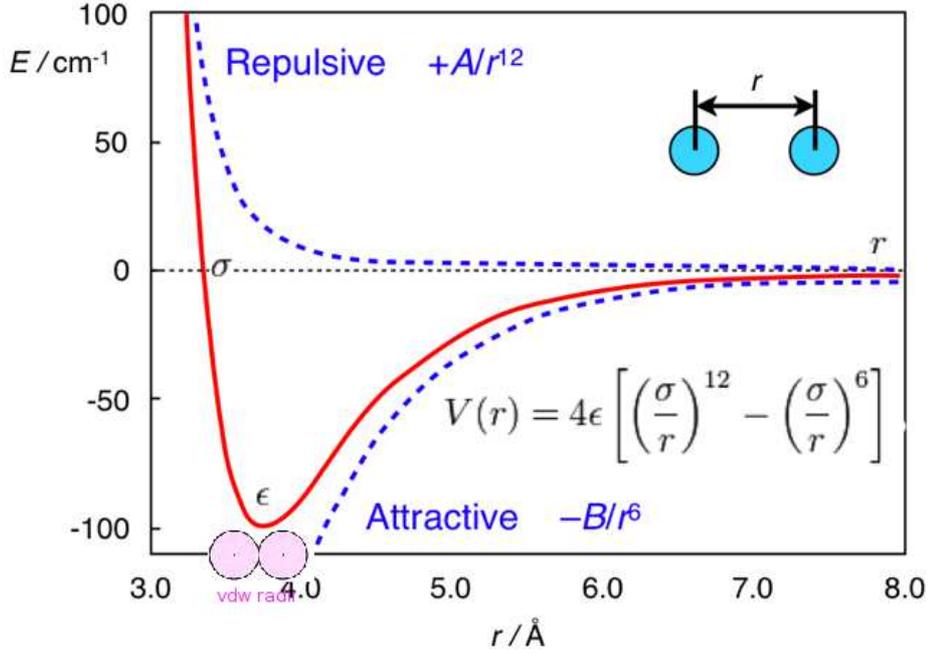}
\caption{{\it The LJ potential energy (Eqs. 1.1 and 1.4) (can be seen in Figure 1 of \cite{zhang_arXiv1106.1584}).}}
\end{figure}
If we introduce the coordinates of the atoms whose number is denoted by $N$ and let $\epsilon =\sigma =1$ be the reduced units, then, for $N$ atoms, Eq. 1.1 is
$$
f(x) = 4\sum_{i=1}^N \sum_{j=1, j<i}^N \left( \frac{1}{\tau_{ij}^6} - \frac{1}{\tau_{ij}^3} \right), \eqno(1.2)
$$ 
where $\tau_{ij} = (x_{3i-2} -x_{3j-2})^2 + (x_{3i-1}-x_{3j-1})^2 +(x_{3i} -x_{3j})^2$, and $(x_{3i-2}, x_{3i-1}, x_{3i})$ is the coordinates of atom $i$, $N\geq 2$. The minimization of LJ potential energy $f(x)$ on $\mathbb{R}^n$ (where $n=3N$) is an optimization problem
$$
\min_{x\in \mathbb{R}^n} f(x), \eqno(1.3)
$$ 
which is a well-known and challenging test problem for global optimization (see {\small http://www-wales.ch.cam.ac.uk/CCD.html} and its recent references such as \cite{kolossváry-and-bowers2010, sicher_etal2010, strodel_etal2010, ye_etal2011, zhang_etal2011b}). 
It is very hard for global optimization to directly solve Eq. 1.3 even without a large number of atoms. Similarly as Eq. 1.1 - i.e.
$$
V_{LJ} (r) = \frac{A}{r^{12}} -\frac{B}{r^6}, \eqno(1.4)
$$
the potential energy for the HBs between $\beta$-strands is
$$
V_{HB} (r) = \frac{C}{r^{12}} -\frac{D}{r^{10}}, \eqno(1.5)
$$
where $A, B, C, D$ are given constants and usually most of the HBs are still kept during the phase of molecular modeling. Thus, the amyloid fibril molecular modeling problem can be reduced to solve the optimization problem Eq. 1.3 though it is not easy to accurately solve Eq. 1.3 for a large molecule.\\

Alternatively, we have found another way to solve Eq. 1.3 \cite{zhang_etal2011a}. Seeing Figure 1, we may know that the optimization problem Eq. 1.3 reaches its optimal value at the bottom of the LJ potential well, where the distance between two atoms equals to the sum of vdW radii of the two atoms. Hence, the amyloid fibril molecular modeling problem can be looked as a molecular distance geometry problem (MDGP) \cite{grosso_etal2009}. As an example to explain MDGP, the problem of locating sensors in telecommunication networks is a DGP. In such a case, the positions of some sensors are known (which are called anchors) and some of the distances between sensors (which can or cannot be anchors) are known. The DGP is to locate the positions of all the sensors. The MDGP looks sensors as atoms and their telecommunication network as a molecule. In mathematics, the following Eqs. 1.6$\sim$1.8 can express the MDGP for Eq. 1.3. The 3D structure of a molecule with $N$ atoms can be described by specifying the 3D coordinate positions $x_1, x_2, \dots, x_N \in \mathbb{R}^3$ of all its atoms. Given bond lengths $d_{ij}$ between a subset $S$ of the atom pairs, the determination of the molecular 3D structure is
$$
(\mathcal{P}_0 ) \quad to \quad find \quad x_1,x_2,\dots ,x_N  \quad such \quad that \quad ||x_i-x_j||=d_{ij}, (i,j)\in S, \eqno(1.6)
$$
where $||\cdot ||$ denotes a norm in a real vector space and in this paper it is calculated as the Euclidean distance 2-norm. Eq. 1.6 can be reformulated as a mathematical global optimization problem
$$
(\mathcal{P} ) \quad \min P(X)=\sum_{(i,j)\in S} w_{ij} (||x_i-x_j||^2 -d_{ij}^2 )^2 \eqno(1.7)
$$
in the terms of finding the global minimum of the function $P(X)$, where $w_{ij}, (i,j)\in S$ are positive weights, $X = (x_1, x_2, \dots, x_N)^T \in \mathbb{R}^n$ \cite{more-and-wu1997} and usually $S$ has fewer elements than $N^2/2$ due to the error in the theoretical or experimental data \cite{zou_etal1997,grosso_etal2009}. Even there may not exist any solution $x_1, x_2, \dots, x_N$ to satisfy the distance constraints in Eq. 1.6, for example when data for atoms $i, j, k \in S$ violate the triangle inequality; in this case, we may add a perturbation item $-\varepsilon^TX$ to $P(X)$:
$$
(\mathcal{P}_{\varepsilon} ) \quad \min P_{\varepsilon}(X)=\sum_{(i,j)\in S} w_{ij} (||x_i-x_j||^2 -d_{ij}^2 )^2 -\varepsilon^TX, \eqno(1.8)
$$
where $\varepsilon \geq 0$. In some cases, instead exact values $d_{ij}, (i,j)\in S$ can be found, we can only specify lower and upper bounds on the distances: $l_{ij} \leq ||x_i - x_j || \leq u_{ij} , (i, j) \in S$; in such cases we may penalize all the unsatisfied constraints into the objective function of ($\mathcal{P}_{\varepsilon}$) by adding
$\sum_{(i,j)\in S} \left( \max \left\{ l^2_{i,j}- ||x_i - x_j ||^2,  0\right\} \right)^2 
                 + \left( \max \left\{ ||x_i - x_j ||^2 - u^2_{i,j}, 0\right\} \right)^2$ into $P_{\varepsilon}(X)$ \cite{zou_etal1997,grosso_etal2009}, where we may let $d_{ij}$ be the interatomic distance (less than 6 angstroms) for the pair in successive residues of a protein and set $l_{ij}=(1-0.05)d_{ij}$ and $u_{ij}=(1+0.05) d_{ij}$ \cite{grosso_etal2009}. In this paper, we aim to solve Eq. 1.8 (or Eq. 1.3) for modeling amyloid fibril molecular 3D structures.\\
                 
Neurodegeneration is the progressive loss of structure or function of neurons, including death of neurons. A prion is a misshapen protein that acts like an infectious agent (but not requiring either DNA, RNA, or both) to cause a number of fatal diseases. Prion diseases are rich in $\beta$-sheets (compared with the normal prion protein PrP$^{\text{C}}$ in rich of $\alpha$-helices) and are so-called ``protein structural conformational" diseases. The normal hydrophobic region 113--120 AGAAAAGA peptide of prion proteins is an inhibitor/blocker of prion diseases. PrP lacking this palindrome could not convert to prion diseases. Brown et al. pointed out that the AGAAAAGA peptide was found to be necessary (though not sufficient) for blocking the toxicity and amyloidogenicity of PrP 106--126, and the peptide AGAA does not form fibrils \cite{brown2000}. The minimum sequence necessary for fibril formation should be AGAAA, AGAAAA, AGAAAAG, AGAAAAGA and GAAAAGA, but the molecular structures of these fibrils have not known yet. This paper addresses an important problem on modeling the 3D molecular structures of prion AGAAAAGA amyloid fibrils of neurodegenerative diseases. The rest of this paper is arranged as follows. In the next section, i.e. Section 2, an improved LBFGS Quasi-Newtonian method is presented for solving Eq. 1.8. Section 3 implements this Quasi-Newtonian method by constructing an 3D molecular structure of prion AGAAAAGA amyloid fibrils of neurodegenerative prion diseases. Numerical results of computations show that the method designed in Section 2 is very effective and successful. This concluding remark will be made in the last section, i.e. Section 4.

\section{Methods}
In a (macro)molecular system, if it is very far from equilibrium, then the forces may be excessively large, a robust energy minimization (EM) is required; another reason to perform an EM is the removal of all kinetic energy from the system: EM reduces the thermal noise in the structures and potential energies \cite{vanderspoel_etal2010}. EM, with the images at the endpoints fixed in space, of the total system energy provides a minimum energy path. EM can be done using steepest descent (SD), conjugate gradient (CG), and Limited-memory Broyden-Fletcher-Goldfarb-Shanno (LBFGS) methods.\\

Three kinds of possible EM methods are: (1) derivative-free methods - that require only function evaluations, e.g. the simplex method and its variants; (2) derivative information methods - the partial derivatives of the potential energy with respect to all coordinates are known and the forces are minimized, e.g. SD, CG methods; and (3) second derivative information methods, e.g. LBFGS method. ``SD is based on the observation that if the real-valued function $f(x)$ is defined and differentiable in a neighborhood of a point $x_0$ then $f(x)$ decreases fastest if one goes from $x_0$ in the direction of the negative gradient of $f(x)$ at $x_0$ and SD local search method  converges fast \cite{snyman2005}. SD is robust and easy to implement but it is not most efficient especially when closer to minimum; at this moment, we may use the efficient CG. CG is slower than SD in the early stages but more efficient when closer to minimum. CG algorithm adds an orthogonal vector to the current direction of the search, and then moves them in another direction nearly perpendicular to this vector. The hybrid of SD-CG will make SD or CG more efficient than SD or CG alone. However, CG cannot be used to find the EM path, for example, when ``forces are truncated according to the tangent direction, making it impossible to define a Lagrangian" \cite{chu_etal2003, case_etal2010}. In this case, the powerful and faster quasi-Newtonian method (e.g. the LBFGS quasi-Newtonian minimizer) can be used \cite{chu_etal2003, liu-and-nocedal1989,nocedal-and-morales2000, byrd_etal1995, byrd_etal19951, zhu_etal1997}. 
We briefly introduce the LBFGS quasi-Newtonian method as follows.\\

Newton's method in optimization explicitly calculates the Hessian matrix of the second-order derivatives of the objective function and the reverse of the Hessian matrix \cite{dennis_etal1996}. The convergence of this method is quadratic, so it is faster than SD or CG. In high dimensions, finding the inverse of the Hessian is very expensive. In some cases, the Hessian is a non-invertible matrix, and furthermore in some cases, the Hessian is symmetric indefinite. Quasi-Newton methods thus appear to overcome all these shortcomings.\\

Quasi-Newton methods (a special case of variable metric methods) are to approximate the Hessian. Currently, the most common quasi-Newton algorithms are the SR1 formula, the BHHH method, the widespread BFGS method and its limited/low-memory extension LBFGS, DFP, MS, and Broyden's methods \cite{nocedal_etal1999, berndt_etal1974, luenberger_etal2008, wang_etal2012}. In Amber \cite{case_etal2010} and Gromacs \cite{vanderspoel_etal2010}, LBFGS is used, and the hybrid of LBFGS with CG - a Truncated Newton linear CG method with optional LBFGS Preconditioning \cite{nocedal-and-morales2000} - is used in Amber \cite{case_etal2010}.\\

\begin{figure}[h!] \label{CPU_time}
\centerline{
\includegraphics[width=4.2in]{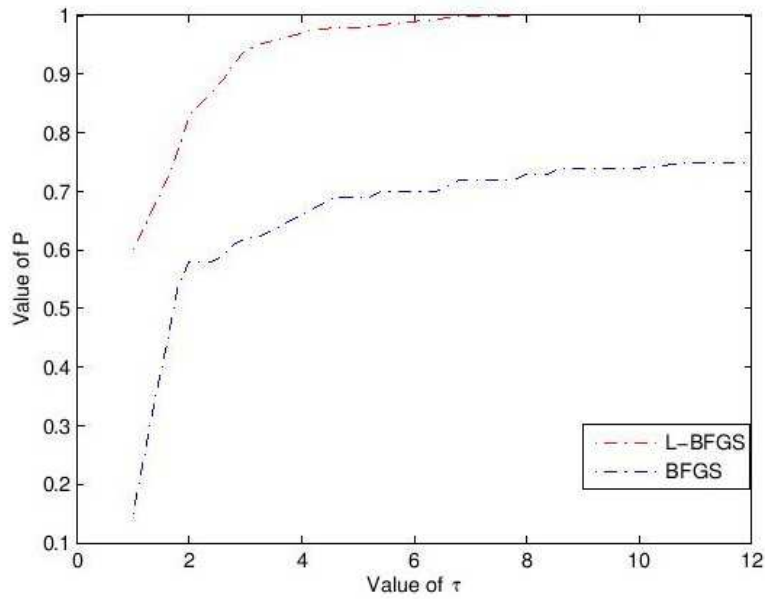}
}
\caption{{\it Performance based on CPU time, where the definitions of P and $\tau$ can be found in \cite{dolan_etal2002}}.}
\end{figure}

For BFGS method, whether it converges at all on nonconvex problems is still an open problem. In fact, Powell (1984) gave a counter-example that shows that BFGS with an inexact line-search search may fail to converge \cite{powell1984, mascarenhas2004, dai2002}. Li and Fukushima (2001) proposed a modified BFGS method for nonconvex objective function \cite{li-and-fukushima2001}. Basing
on \cite{li-and-fukushima2001, liu-and-nocedal1989, xiao_etal2008, nocedal-and-morales2000, yang-and-xu2007}, in this paper we present an improved LBFGS method described as follows \cite{hou_etal2012} - which presents the nonmonotone line search technique \cite{grippo_etal1986, han_etal1997} for the Wolfe-type search. The improved LBFGS method presented in this paper is much better than the standard BFGS method in view of the CPU time (see Figure \ref{CPU_time}) tested through more than 30 nonlinear programming problems (where each selected problem is regular, that is, its first and second derivatives exist and are continuous everywhere, and each problem is with different dimensions, i.e., 100, 500, 1000 and 10000 dimensions) and its mathematical theory to support this algorithm can be seen from the Supplementary Materials \cite{hou_etal2012} listed at the end of this paper. This paper implements the Wolfe-type search by the approximation technique of piecewise linear/quadratic function \cite{bagirov_etal2008}.

\begin{algorithm} {\bf An Improved LBFGS Method for minimizing nonconvex function \cite{hou_etal2012}}
\item[] {\it Step 0:} Choose an initial point $x_0\in \mathbb{R}^n$, an initial positive definite matrix $H_0$, and choose constants $\sigma_1$,$\sigma_2$ such that $0<\sigma_1<\sigma_2<1$, and choose an positive integer $m_1$. Let $k=0$.
\item[] {\it Step 1:} If $\|g_k\|=0$, then output $x_k$ and stop; otherwise, go to {\it Step 2}.
\item[] {\it Step 2:} Solve the following linear equation to get $d_k$:
$$d_k=-H_kg_k,$$
\item[] {\it Step 3:} Find a step-size $\lambda_k>0$ satisfying the Wolfe-type line search conditions:
$$f(x_k+\lambda_kd_k)\leq f(x_k)+\sigma_1\lambda_kg_{k}^Td_k, \eqno(2.1)$$
$$g(x_k+\lambda_kd_k)^Td_k\geq\sigma_2 g_{k}^Td_k. \eqno(2.2)$$
\item[] {\it Step 4:} Let the next iterate by $x_{k+1}=x_k+\lambda_kd_k$. Calculate $g_{k+1}$ and $\| g_{k+1} \|$.
\item[] {\it Step 5:} Let $s_k=x_{k+1}-x_k=\lambda_kd_k$,$y_k=g_{k+1}-g_k$, $\gamma_k=\|g_k\|$, then $y_{k}^\ast=y_k+\gamma_ks_k$.
\item[] {\it Step 6:} Let $\bar {m}=\min\{k+1,m_1\}$. Update $H_k$ following the formula
$$\begin{array}{rcl}
H_{k+1}&=&(V_{k}^{\ast T}\cdots V_{k-\bar {m}+1}^{\ast T})H_{k-\bar {m}+1}(V_{k-\bar {m}+1}^\ast\cdots V_{k}^\ast)\\
&&+\omega_{k-\bar {m}+1}^\ast(V_{k-1}^{\ast T}\cdots V_{k-\bar {m}+2}^{\ast T}) s_{k-\bar {m}+1}s_{k-\bar {m}+1}^T(V_{k-\bar {m}+2}^\ast\cdots V_{k-1}^\ast)\\
&&+\cdots+\omega_{k}^\ast s_ks_{k}^T,
\end{array} \eqno(2.3)$$
where $\omega_{k}^\ast=\frac{1}{s_{k}^Ty_{k}^\ast}$, $V_{k}^\ast=I-\omega_{k}^\ast y_{k}^\ast s_{k}^T$ (when $k=0$, $H_0=\frac{y_0^Ts_0}{\|y_0\|^2}, V_0^*=I-\frac{y_0^*s_0^T}{s_0^Ty_0^*}, H_1=(I-\frac{y_0^*s_0^T}{s_0^Ty_0^*})^TH_0(I-\frac{y_0^*s_0^T}{s_0^Ty_0^*}) +\frac{1}{s_0^Ty_0^*} s_ks_k^T$). 
\item[] {\it Step 7:} Let $k=k+1$ and go to {\it Step 1}.
\end{algorithm}
\noindent In Algorithm 1, $g_k$ denotes the gradient of $f(x)$ at $x_k$, and the convergence and the R-linear convergent rate of Algorithm 1 are guaranteed by the following assumptions: (A1) $f$ is twice continuously differentiable, (A2) $f$ has Lipschitz continuous gradients and Hessians, (A3) the Hessian at the stationary point is always positively definite, and (A4) the level set of $f$ is bounded. The detailed proof of the convergence and $R$-linear convergent rate of Algorithm 1 can be found in the Supplementary Materials supplied at the end of this paper.

In \cite{hou_etal2012}, the following nonmonotone line search technique for stepsize $\lambda_k$ is used:
$$f(x_k+\lambda_kd_k)\leq \max_{0\leq j\leq M_0} f(x_{k-j})+\sigma_1\lambda_kg_k^Td_k,$$
$$g(x_k+\lambda_kd_k)^Td_k \geq \max \{ \sigma_2, 1-(\lambda_k ||d_k||)^p\}  g_{k}^Td_k,$$
where $p\in (-\infty , 1)$, and $M_0$ is a nonnegative integer. This is a difference between the algorithm of \cite{hou_etal2012} and this paper. All in all, it is well known that quasi-Newton method is an efficient solution method for unconstrained and continuously differentiable minimization problem \cite{byrd_etal1989, byrd_etal1987, dennis_etal1977}. However, it needs computing and storage of the updated matrix which is an approximation to the Hessian matrix at each iteration of the method. Hence, its efficiency may decrease when it is applied to large scale optimization problem. To overcome the drawback, limited memory quasi-Newton method is proposed \cite{nocedal1980}. The main ideal of this method is nearly identical to that of the standard BFGS method, and the only difference is that the inverse Hessian approximation is not formed explicitly, but defined by a small number, say $\bar{m}$, of BFGS updates. This technique received much attention in recent years and numerical experiments show that it is very competitive \cite{gilbert_etal1989, liu-and-nocedal1989}, and its global convergence and R-linear convergence rate with Wolfe line search are established for the uniformly convex case \cite{byrd_etal1995, liu-and-nocedal1989}. Since the limited memory BFGS method may suffer from ill-conditions for small value of $\bar{m}$, Al-Baali (1999) \cite{albaali1999} made some modifications to the method and establish its global convergence based on the same assumptions, and Byrd et al. (1994) \cite{byrd_etal1994} derives new representation of limited memory quasi-Newton matrices for the benefit of computing the updated matrix. Recently, a non-monotone line search is introduced, see e.g., \cite{grippo_etal1986, han_etal1997}. Then it is showed to be more competitive and practical for solving nonlinear optimization problems, and \cite{yuan_etal2010} established the global convergence of this line search applied to limited memory BFGS method based on the uniformly convex assumption. Motivated by the above observation, it turn out that in two respects the limited memory BFGS method is much less effective. First, we note that the convergence analysis of these method are focused on the uniformly convex assumption and little is known for nonconvex case. Second, numerical experiments have suggested the main weakness of limited memory method is that it may converge very slowly in terms of number of iterations for ill-conditioned problems. The purpose of the above Algorithm 1 is to reduce these defects and Figure 2 shows the effectiveness of the proposed algorithm. We will apply it into the molecular modeling of prion AGAAAAGA amyloid fibrils in the next section.

\section{Results and Discussion}
From their research of prion, scientists found that the cross-$\beta$ structure of peptides is with the nature of self-aggregation, the self-aggregating to form fibers. This provides us a new research idea for nanomaterials. HBs can be formed between peptide $\beta$-strands, and one peptide monomer connects together with another in accordance with the specific structure to form fibers. Many laboratories in the world are synthesizing peptides that can self-aggregate to form fibers, and want to be able to control the growth of the fiber to find out new functional materials \cite{alper_etal1967, griffith1967}. The studies of this paper not only benefit nanometerials research, but also benefit the research on neurodegenerative amyloid fibril diseases. Prion AGAAAAGA peptide has been reported to own an amyloid fibril forming property (initially described in 1992 by Gasset et al. of Prusiner's Group)
\cite{brown2000, brown2001, brown_etal1994, cappai-and-collins2004, chabry_etal1998, cheng_etal2003, gasset_etal1992, haigh_etal2005, harrison_etal2010, holscher_etal1998, jobling_etal2001, jobling_etal1999, jones_etal2011, kourie2001, kourie_etal2003, kuwata_etal2003, laganowsky_etal2012, lee_etal2008, ma_etal2002, norstrom-and-mastrianni2005, sasaki_etal2008, wagoner2010, wagoner_etal2011, wegner_etal2002, zanuy_etal2003}, 
but there has not been experimental structural bioinformatics for this segment yet due to the unstable, noncrystalline and insoluble nature of this region. Furthermore, Zhang (2011) did accurate calculations to confirm the amyloid fibril property at this region (Figure 3) \cite{zhang2011}.\\

\begin{figure}[h!] \label{Fig03-identification}
\includegraphics[width=4.2in]{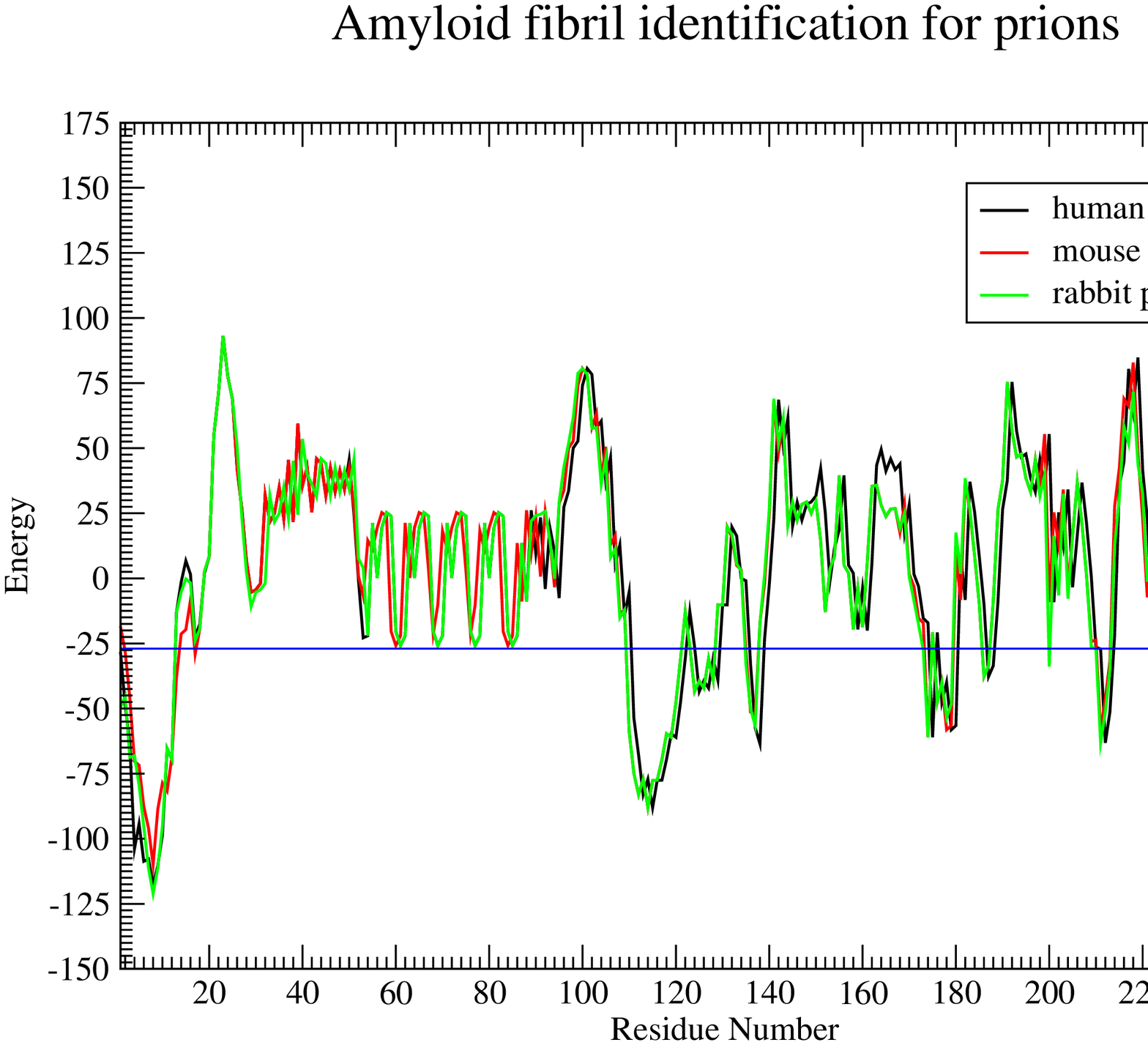}
\caption{{\it Prion AGAAAAGA (113-120) segment is clearly and surely identified as the amyloid fibril formation region, because its energy is less than the amyloid fibril formation threshold energy -26 kcal/mol \cite{zhang_etal2007}.}}
\end{figure}     

\subsection{Material for the Molecular Modeling}
This paper uses a suitable pdb file template 2OMP.pdb (the LYQLEN peptide derived from human insulin residues 13-18 \cite{sawaya_etal2007}) from the Protein Data Bank to build an 8-chain AGAAAAGA prion amyloid fibril molecular model to illuminate Algorithm 1 works very well. To choose 2OMP.pdb (Figure 4) as the modeling template is due to it can pass all the long procedures of SDCG-SA (equilibrations \& productions)-SDCG of \cite{zhang2011}.
\begin{figure}[h!] \label{Fig04-2OMP-pdb}
\includegraphics[width=5.5in]{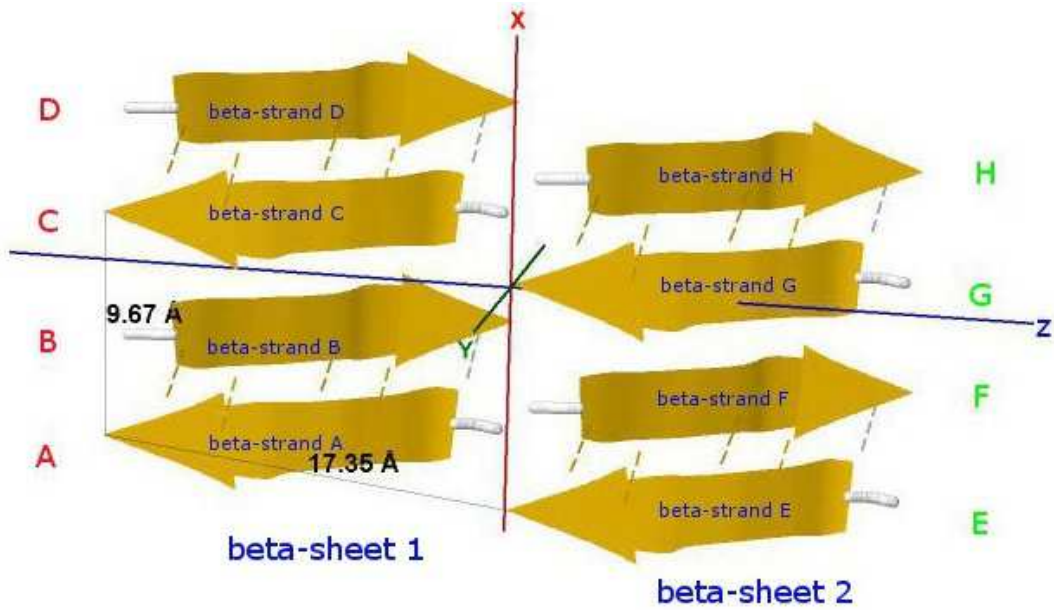}
\caption{{\it Protein fibril structure of human insulin LYQLEN (13-18) (PDB id: 2OMP). The dashed lines denote the HBs between the pairs of $\beta$-strands. A, B, C, D, E, F, G, H denote the chains of the fibril. The pair of $\beta$-sheets 1 \& 2 forms a completely dry interface by vdWs, and between many pairs of $\beta$-sheets wet interfaces are formed with water molecules.}}
\end{figure}     
By observations of Figure 4 and the 2nd column of coordinates of 2OMP.pdb, we know that E(F) chains can be calculated on the XZ-plane from A(B) chains by Eq. 3.1 and other chains can be got by a parallel up (or down) along the X-axis by Eqs. 3.2$\sim$3.3:
$$E(F) = A(B) + (-1.885, 0, 17.243), \eqno(3.1)$$
$$C(D) = A(B) + ( 9.666,0,0), \eqno(3.2)$$
$$G(H) = E(F) + ( 9.666,0,0). \eqno(3.3)$$

\subsection{New Molecular Modeling Homology Model}
Basing on the template 2OMP.pdb from the Protein Data Bank (Figure 4), the AGAAAAGA palindrome amyloid fibril model of prions (denoted as Model 1) will be constructed. Chains AB of Model 1 will be got from AB Chains of 2OMP.pdb using the mutate module of the free package Swiss-PdbViewer (SPDBV Version 4.01) (http://spdbv.vital-it.ch). It is pleasant to see that some HBs are still kept after the mutations; thus we just need to consider the vdWs only. Making mutations for EF Chains of 2OMP.pdb, we can get the EF Chains of Model 1. Then we add GLY and ALA residues by XLEaP module of Amber 11. However, the vdWs between Chain A and Chain E, between B Chain and F Chain are too far at this moment (Figure 5, where the twice of the vdW radius of CB atom is 3.4 angstroms).\\ 

\begin{figure}[h!] \label{Fig05-AGAAAAGA-ABEF}
\includegraphics[width=3.902in]{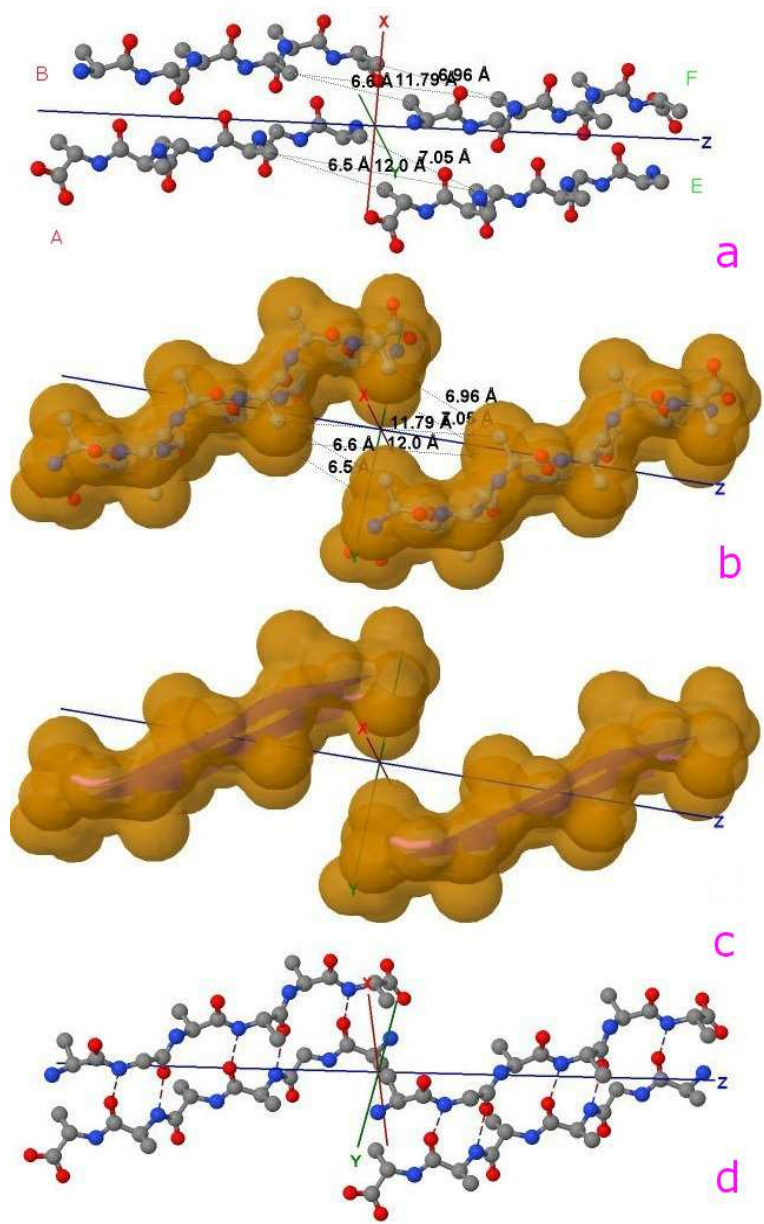}
\caption{{\it 5a shows the distances of ``Zipper 1"-E.ALA6.CB-A.ALA3.CB-E.ALA4.CB-A.ALA1.CB are 6.5, 12.0 and 7.05 angstroms respectively, and the distances of ``Zipper 2"-F.ALA1.CB-B.ALA4.CB-F.ALA3.CB-B.ALA6.CB are 6.6, 11.79, 6.96 angstroms respectively. 5b shows the far vdW surface. 5c shows the violet colored ABEF Chains of Figure 3. 5d shows HBs: A/E.ALA5.O-B/F.GLY2.N, A/E.ALA5.N-B/F.GLY2.O, A/E.ALA3.O-B/F.ALA4.N, A/E.ALA3.N-B/F.ALA4.O, A/E.ALA1.O-B/F.ALA6.N.}}
\end{figure}

In \cite{zhang2011} the commercial package InsightII (http://accelrys.com/) is used to build models. Instead of InsightII, because this package is not available by the authors, this paper uses Algorithm 1 to build and optimize Model 1. In ``Zipper 1", fixing the coordinates of A.ALA3.CB, A.ALA1.CB and letting the coordinates of E.ALA6.CB, E.ALA4.CB be variables, we may get an optimization problem:
{\small
$$\begin{array}{rcl}
\min 4 \left\{ \frac{1}{ \left[ (x_{11}-1.071)^2 +(x_{12}-2.986)^2 +(x_{13}-1.888)^2 \right]^6 } - 
               \frac{1}{ \left[ (x_{11}-1.071)^2 +(x_{12}-2.986)^2 +(x_{13}-1.888)^2 \right]^3 } \right\} \\
    +4 \left\{ \frac{1}{ \left[ (x_{21}-1.071)^2 +(x_{22}-2.986)^2 +(x_{23}-1.888)^2 \right]^6 } - 
               \frac{1}{ \left[ (x_{21}-1.071)^2 +(x_{22}-2.986)^2 +(x_{23}-1.888)^2 \right]^3 } \right\}\\
    +4 \left\{ \frac{1}{ \left[ (x_{21}-1.135)^2 +(x_{22}+0.763)^2 +(x_{23}-7.209)^2 \right]^6 } - 
               \frac{1}{ \left[ (x_{21}-1.135)^2 +(x_{22}+0.763)^2 +(x_{23}-7.209)^2 \right]^3 } \right\}           
\end{array} \eqno(3.4)$$ 
}
or
{\small
$$\begin{array}{rcl}
\min \frac{1}{2} \left\{  (x_{11}-1.071)^2 +(x_{12}-2.986)^2 +(x_{13}-1.888)^2 -3.4^2 \right\}^2\\
    +\frac{1}{2} \left\{  (x_{21}-1.071)^2 +(x_{22}-2.986)^2 +(x_{23}-1.888)^2 -3.4^2 \right\}^2\\
    +\frac{1}{2} \left\{  (x_{21}-1.135)^2 +(x_{22}+0.763)^2 +(x_{23}-7.209)^2 -3.4^2 \right\}^2\\
    -0.05 \left\{ x_{11} +x_{12} +x_{13} +x_{21} +x_{22} +x_{23} \right\}            
\end{array} \eqno(3.5)$$
}
with an initial solution (-0.067, 5.274, 7.860; -1.119, 1.311, 13.564). Similarly, in ``Zipper 2", fixing the coordinates of B.ALA4.CB, B.ALA6.CB and letting the coordinates of F.ALA1.CB, F.ALA3.CB be variables, we get another optimization problem:
{\small
$$\begin{array}{rcl}
\min 4 \left\{ \frac{1}{ \left[ (x_{11}-5.446)^2 +(x_{12}-2.796)^2 +(x_{13}-2.662)^2 \right]^6 } - 
               \frac{1}{ \left[ (x_{11}-5.446)^2 +(x_{12}-2.796)^2 +(x_{13}-2.662)^2 \right]^3 } \right\} \\
    +4 \left\{ \frac{1}{ \left[ (x_{21}-5.446)^2 +(x_{22}-2.796)^2 +(x_{23}-2.662)^2 \right]^6 } - 
               \frac{1}{ \left[ (x_{21}-5.446)^2 +(x_{22}-2.796)^2 +(x_{23}-2.662)^2 \right]^3 } \right\}\\
    +4 \left\{ \frac{1}{ \left[ (x_{21}-5.201)^2 +(x_{22}+1.125)^2 +(x_{23}-7.873)^2 \right]^6 } - 
               \frac{1}{ \left[ (x_{21}-5.201)^2 +(x_{22}+1.125)^2 +(x_{23}-7.873)^2 \right]^3 } \right\}           
\end{array} \eqno(3.6)$$ 
}
or
{\small
$$\begin{array}{rcl}
\min \frac{1}{2} \left\{ (x_{11}-5.446)^2 +(x_{12}-2.796)^2 +(x_{13}-2.662)^2 -3.4^2 \right\}^2\\
    +\frac{1}{2} \left\{ (x_{21}-5.446)^2 +(x_{22}-2.796)^2 +(x_{23}-2.662)^2 -3.4^2 \right\}^2\\
    +\frac{1}{2} \left\{ (x_{21}-5.201)^2 +(x_{22}+1.125)^2 +(x_{23}-7.873)^2 -3.4^2 \right\}^2\\
    -0.05 \left\{ x_{11} +x_{12} +x_{13} +x_{21} +x_{22} +x_{23} \right\}            
\end{array} \eqno(3.7)$$
}
with an initial solution (4.714, 4.878, 8.881; 4.170, 1.360, 14.292). Next, we solve Eqs. 3.5 and 3.7 by Algorithm 1.\\

We first solve Eq. 3.5 in the use of Algorithm 1. We set $\sigma_1=10^{-4}$, $\sigma_2=0.1\sim 0.9$, $m_1=3\sim 7$, take the initial solution $x_0=$(-0.067, 5.274, 7.860; -1.119, 1.311, 13.564) and calculate its gradient $g_0$=(-69.7747, 140.135, 365.852, -752.075, -285.005, 3576.69) and its Hessian matrix $H_0=$
{\small
$$
\left( \begin{array}{cccccc}
66.4497  &-10.415  &-27.1845 &0           &0            &0\\
-10.415  &82.2093  &54.6557  &0           &0            &0\\
-27.1845 &54.6557  &203.929  &0           &0            &0\\
0        &0        &0        &380.664     &-4.02618     &-159.578\\
0        &0        &0        &-4.02618    &369.586      &-25.5081\\
0        &0        &0        &-159.578    &-25.5081     &1048.02 \end{array} \right),
\eqno(3.8)$$
}
which is a positive definite matrix with eigenvalues (1085.02, 372.093, 341.157, 230.049, 61.2695, 61.2695). Then Algorithm 1 hybridized with simulated annealing global optimal search (in order to bring local optimal solutions to jump out of local traps, replacing the discrete gradient local optimal search method in Algorithm 1 of \cite{zhang_etal2011b} by the Algorithm 1 of this paper) is executed and the optimal solution (3.027, 4.954, 3.856; 1.679, 1.777, 5.011) for Eq. 3.5 is got.\\

Similarly, for Eq. 3.7, we take the initial solution $x_0=$(4.714, 4.878, 8.881; 4.170, 1.360, 14.292) and calculate its gradient $g_0$=(-46.8782, 133.142, 397.798, -401.192, -182.604, 3436.46) and its Hessian matrix $H_0=$
{\small
$$
\left( \begin{array}{cccccc}
66.1163  &-6.0961  &-18.2092     &0           &0            &0\\
-6.0961  &81.3119  &51.7918      &0           &0            &0\\
-18.2092 &51.7918  &218.677      &0           &0            &0\\
0        &0        &0            &339.302     &-2.9188      &-85.8315\\
0        &0        &0            &-2.9188     &361.487      &-2.99786\\
0        &0        &0            &-85.8315    &-2.99786     &1034.38 \end{array} \right),
\eqno(3.9)$$
}
which is a positive definite matrix with eigenvalues (1044.83, 361.8, 328.538, 238.159, 63.973, 63.973). The optimal solution (7.412, 4.760, 4.624; 5.887,  1.451, 5.757) for Eq. 3.7 is got.\\
 
We set (3.027, 4.954, 3.856; 1.679, 1.777, 5.011) as the coordinates of E.ALA6.CB and E.ALA4.CB, (7.412, 4.760, 4.624; 5.887,  1.451, 5.757) as the coordinates of F.ALA1.CB and F.ALA3.CB, and taking the average value we get
$$E(F) = A(B) + (1.0335, 1.0823, 0.9723). \eqno(3.10)$$
By Eq. 3.10 we can get very close vdW contacts between A (B) chains and E(F) chains (Figure 5(b)), and other chains of Model 1 can be got by Eqs. 3.2$\sim$3.3 and Eq. 3.11:
$$I(J) = A(E) - ( 9.666,0,0). \eqno(3.11)$$\\

\begin{figure}[h!] \label{Fig6-AGAAAAGA}
\includegraphics[width=5.8in]{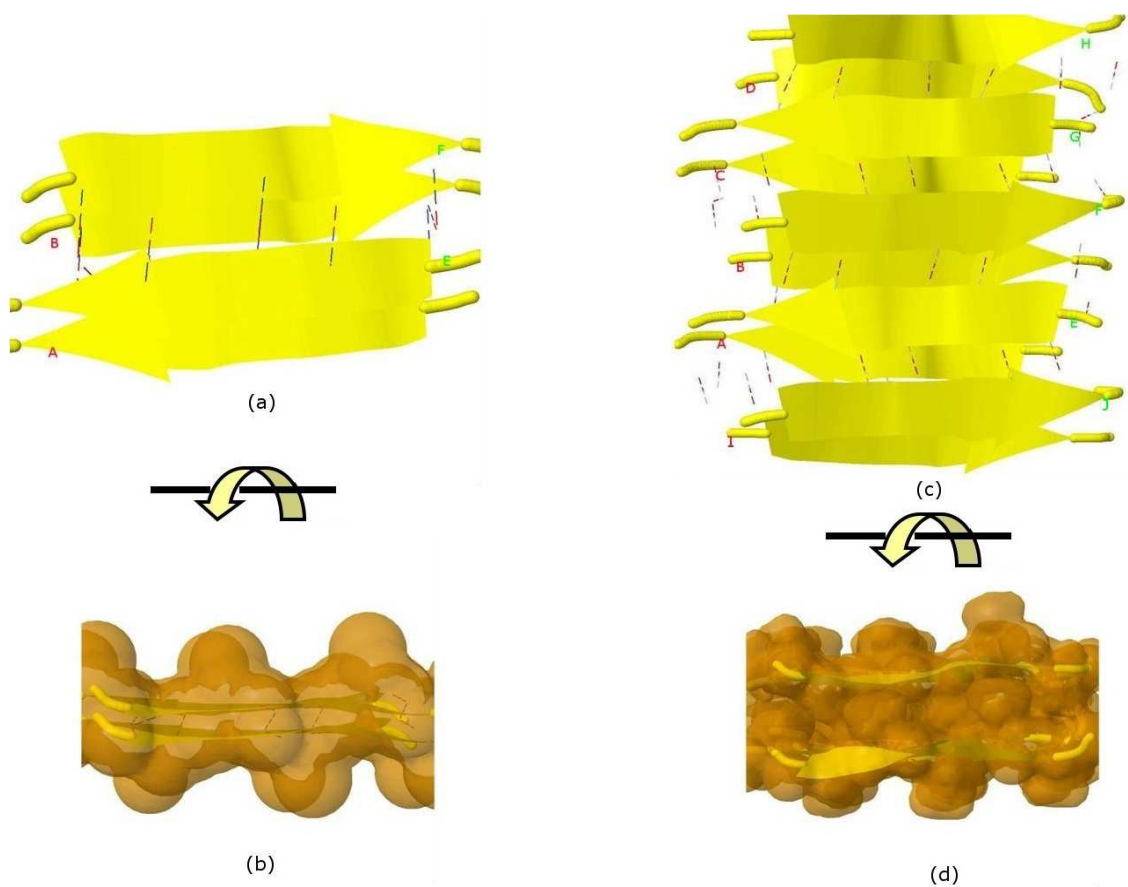} 
\caption{{\it Left: The close vdW surface contacts between Chains A(B) and E(F) of Model 1 after solving Eqs. 3.5 and 3.7; Right: The constructed Model 1 with 10 chains. The dashed lines denote HBs.}}
\end{figure}

The initial structure of Model 1 illuminated in Figure 6(a)$\sim$6(b) – is not the optimal structure with the lowest total potential energy. The initial structure also has no hydrogen atoms (so no hydrogen bonds existed) and water molecules added. For each Chain, the C-terminal and N-terminal atoms also have problems. Clearly there are a lot of close/bad contacts between $\beta$-strand atoms as illuminated in Figure 6(a)$\sim$6(b). We used the ff03 force field of AMBER 11,in a neutral pH environment. The amyloid fibrils were surrounded with a 8 angstroms layer of TIP3PBOX water molecules using the XLEaP module of AMBER 11. 1944 waters and 408 hydrogen atoms were added for Model 1 byt he XLEaP module. solvated amyloid fibril was inimized by the method.

The LJ potential energy of atoms' vdW interactions is just a part of the total potential energy of a protein, and by observations from Model 1 computed by Eqs. 3.10, 3.2$\sim$3.3 and 3.11 we can see that the contacts between $\beta$-strand atoms and $\beta$-sheet atoms are too close/bad. Thus, we need to relax Model 1 computed. The relaxation is done in the use of local search LBFGS Quasi-Newton method (lbfgs\_memory\_depth=3) within AMBER 11 \cite{case_etal2010}. The relaxed/optimized Model 1 is illuminated in Figure 6(c). Seeing Figure 6(d) compared with Figure 6(b), we may know the vdW interactions between the two $\beta$-sheets are very relaxed/optimized now. Figure 6(c) shows the Model 1 of optimal molecular structure for prion AGAAAAGA amyloid fibrils.

\section{Conclusion}
In a (macro)molecular system, a robust energy minimization is very necessarily required. Mathematical optimization minimization methods find a place to apply in these systems. Because in physics the (macro)molecular system usually is not a simple two-body problem of system, local search optimization methods are very useful in the applications to the (macro)molecular system. On anther sense, when a protein is unstable, noncrystalline or insoluble and very difficult to detect its 3D structure by the expensive and costly NMR and X-ray, theoretical mathematical or physical computational method can be used to produce the 3D structure of the protein. Moreover, even the X-ray crystallography finds the X-ray final structure of a protein, we still need refinements using theoretical protocols in order to produce a better structure. The theoretical computational method - an improved LBFGS Quasi-Newtonian mathematical optimization method - presented in this paper and other mathematical optimization methods mentioned in this paper should be very useful in the protein molecular modeling research field.\\

This paper also shows the effectiveness of the improved LBFGS mathematical optimization method presented. Prion AGAAAAGA amyloid fibrils have not much structural information. This paper presents some bioinformatics on the molecular structures of prion AGAAAAGA amyloid fibrils in the sense of theoretical emphasis. The structures may be helpful in the advance in the biochemical knowledge of prion protein misfolding or instability and in the future applications for therapeutic agent design.\\ 

\section*{Acknowledgments} 
This research was supported by the National Natural Science Foundation of China (No. 11171180), the National Basic Research Program of China (2010CB732501), and a Victorian Life Sciences Computation Initiative (VLSCI) grant numbered VR0063 \& 488 on its Peak Computing Facility at the University of Melbourne, an initiative of the Victorian Government (Australia). The first author appreciates the financial support by Professors Changyu Wang and Yiju Wang for his visits to Qufu Normal University (China), and by Professor Xiangsun Zhang for his visits to Chinese Academy of Sciences (Beijing); the first author also appreciates Professor Adil M. Bagirov (University of Ballarat) for his instructions to implement the Wolfe-type search.  This paper is dedicated to Professor Xiangsun Zhang (Academia Sinica, Beijing) on the occasion of his 70th birthday.

\section*{Supplementary Materials \cite{hou_etal2012}}
\subsection*{Mathematical Assumptions (A1)-(A4):}

\vskip 0.5cm

\qquad (A1) the objective function $f:\mathbb{R}^n\rightarrow \mathbb{R}$ is twice continuously differentiable, i.e. $f\in C^2$;

(A2) $f$ has Lipschitz continuous gradients and Hessians, i.e., there exist constants $L, M_2>0$ such that
$$\|g(x)-g(y)\|\leq L\|x-y\|,\quad \forall x,y\in \mathbb{R}^n,$$
$$\|G(x)-G(x^\ast)\|\leq M_2\|x-x^\ast\|,$$
where $G$ denotes the Hessian and $x^\ast$ is a stationary point (i.e. ${x_k}$ converges to $x^\ast$ where $g(x^\ast)=0$);

(A3) $G(x^\ast)$ is positive definite; and

(A4) the level set of $f$
$$\Omega=\{x\in \mathbb{R}^n~|~f(x)\leq f(x_0)\},$$
is bounded, where $x_0\in \mathbb{R}^n$ is an initial point.

\subsection*{Mathematical Proof of the Algorithm Convergence and Convergent Rate:}

\vskip 0.5cm

In order to establish the convergence for Algorithm 1, we represent Eq. 2.3 as
$$B_{k}^{l+1}=B_{k}^l-\frac{B_{k}^ls_ls_{l}^TB_{k}^l}{s_{l}^TB_{k}^ls_l}+\frac{y_{l}^\ast y_{l}^{\ast T}}{y_{k}^{\ast T}s_l},\quad l=k-\bar {m}+1,\cdots,k,\eqno(5.1)$$
where $B_k=H_{k}^{-1}$,$s_l=x_{l+1}-x_l$, $y_{l}^\ast=y_l+\gamma_ls_l$, and $B_{k}^{k-\bar {m}+1}=B_0$ for all $k$.\\

\noindent \begin{theorem}\label{thm51}\quad Let Assumptions (A1)$\sim$(A4) hold, and $\{x_k\}$ denotes the sequence generated by Algorithm 1. Then $\lim\limits_{k\rightarrow\infty}\inf\|g_k\|=0$; moreover, there exits a constant $t\in [0,1)$ such that
$$f(x_k)-f(x^\ast)\leq t^k[f(x_0-f(x^\ast)],$$
which is to say that ${x_k}$ converges to the minimum $R$-linearly.
\end{theorem}
{\it Proof.}\quad Denote by $M_1$ an upper bound of ${\|g_k\|}$ on $\Omega$ and by $y_{k}^\ast=y_k+\gamma_ks_k$, we can prove that
$$\frac{{\|y_{k}^{\ast T}\|}^2}{y_{k}^{\ast T} s_k}\leq\frac{(L+M_1)^2}{\gamma},$$ where $\gamma$ is a sufficiently small positive number. By this formula and taking the trace operation in both sides of Eq. 5.1, we get
$$\begin{array}{rcl}
tr(B_{k+1})&=&tr(B_{k}^{k-\bar {m}+1})-\sum\limits_{l=k-\bar {m}+1}^k\frac{\|B_ls_l\|^2}{s_{l}^TB_ls_l}+\sum\limits_{l=k-\bar {m}+1}^k\frac{{\|y_{l}^\ast\|}^2}{y_{l}^{\ast T}s_l}\\
&=&\cdots\\
&\leq&tr(B_0)-\sum\limits_{l=0}^k\frac{\|B_ls_l\|^2}{s_{l}^TB_ls_l}+\sum\limits_{l=0}^k\frac{{\|y_{l}^\ast\|}^2}{y_{l}^{\ast T}s_l}\\
&\leq&tr(B_0)-\sum\limits_{l=0}^k\frac{\|B_ls_l\|^2}{s_{l}^TB_ls_l}+\frac{(L+M_1)^2}{\gamma}(k+1)\\
&\leq&M_3,
\end{array}\eqno(5.2)$$
where $M_3$ is a positive constant. We also take the determinant and by Eq. 5.2, we have
$$\begin{array}{rcl}
det(B_{k+1})&=&det(B_{k}^{k-\bar {m}+1})\prod\limits_{l=k-\bar {m}+1}^k\frac{y_{l}^{\ast T}}{s_{l}^TB_ls_l}\\
&=&\cdots=det(B_0)\prod\limits_{l=0}^k\frac{y_{l}^{\ast T}}{s_{l}^TB_ls_l}
\geq det(B_0)(\frac{\gamma}{M_3})^{k+1}\\
&\geq&M_4,
\end{array}\eqno(5.3)$$where $M_4$ is a positive constant. Combining Eqs. 5.2 and 5.3, we obtain a constant $\delta>0$ such that $\cos(\theta_k)\geq\delta$, where $\cos(\theta_k)=\frac{s_{k}^TB_ks_k}{\|s_k\|\|B_ks_k\|}$.\\

Because $g(x)$ is Lipschitz continuous and by the Wolfe-type line search condition Eq. 2.2, we have
$$(\sigma_2-1)g_{k}^Td_k\leq(g_{k+1}-g_k)^Td_k\leq\lambda_kL\|d_k\|^2,$$
which implies
$$\lambda_k\geq\frac{(\sigma_2-1)}{L}\frac{g_{k}^Td_k}{\|d_k\|^2}.\eqno(5.4)$$\\

Using the Wolfe-type line search condition Eq. 2.1 and Eq. 5.4, we have
$$\begin{array}{rcl}
f_{k+1}&\leq&f_k+\sigma_1\lambda_kg_{k}^Td_k\leq f_k+\sigma_1\frac{(\sigma_2-1)}{L}\frac{(g_{k}^Td_k)^2}{\|d_k\|^2}\\
&=&f_k+\sigma_1\frac{(\sigma_2-1)}{L}\|g_{k}\|^2\cos^2\theta_k.\end{array}$$
Then using the monotonicity and the bound of $\{f_k\}$ on $\Omega$, we obtain
$$\sum\limits_{k=1}^\infty\|g_k\|^2\cos^2\theta_k<\infty.$$
Combining $\cos\theta_k\geq\delta$, we get
$$\lim_{k\rightarrow\infty}\inf\|g_k\|=0.$$\\

Under the Assumptions (A2)$\sim$(A3), there is a neighbourhood $U(x^\ast)$ of $x^\ast$ such that for all $x\in U(x^\ast)$,
$$\|g(x)\|\geq\|g(x)-g(x^\ast)\|\geq m\|x-x^\ast\|,\eqno(5.5)$$
and for all $d\in R^n$,
$$d^TG(x)d\geq m\|d\|^2,\eqno(5.6)$$
where $m$ is a positive constant got according to the Assumptions. Hence, Eqs. 5.5$\sim$5.6 hold with $x=x_k$ for all $k$ sufficiently large.\\

Since
$$\begin{array}{rcl}
y_{k}^{\ast T}s_k&=&(y_k+\gamma_ks_k)^Ts_k\geq y_{k}^Ts_k\\
&\geq&-(1-\sigma_2)g_{k}^Ts_k\\
&=&(1-\sigma_2)\|g_k\|\|s_k\|\cos\theta_k,\end{array}$$
and $y_{k}^{\ast T}s_k\leq(L+M_1)^2\|s_k\|^2$, if denoting $\alpha_1=\frac{1-\sigma_1}{(L+M_1)^2}$, we can obtain
$$\|s_k\|\geq\frac{1-\sigma_1}{(L+M_1)^2}\|g_k\|\cos\theta_k=\alpha_1\|g_k\|\cos\theta_k.\eqno(5.7)$$
Combining Eq. 5.7 and the Wolfe-type line search condition Eq. 2.1, we get
$$\begin{array}{rcl}
f_{k+1}-f^\ast&\leq&f_k-f^\ast+\sigma_1g_{k}^Ts_k\\
&=&f_k-f^\ast-\sigma_1\|g_k\|\|s_k\|\cos\theta_k\\
&\leq&f_k-f^\ast-\sigma_1\alpha_1\|g_k\|^2\cos^2\theta_k.\end{array}$$
Using Eq. 5.5,  then Eq. 5.7 implies for all $k$ sufficiently large,
$$\begin{array}{rcl}
f_{k+1}-f^\ast&\leq&f_k-f^\ast-\sigma_1\alpha_1\|g_k\|^2\cos^2\theta_k\\
&\leq&f_k-f^\ast-\sigma_1\alpha_1m^2\|x_k-x^\ast\|^2\cos^2\theta_k.
\end{array}\eqno(5.8)$$
By Taylor's expansion, then there exists a positive number $M_5>0$ such that
$$f_k-f^\ast\leq M_5\|x_k-x^\ast\|^2$$
for all $k$ sufficiently large.\\

Let $\alpha_2=\frac{\sigma_1\alpha_1m^2}{M_5}$, indeed Eq. 5.8 implies that
$$\begin{array}{rcl}
f_{k+1}-f^\ast&\leq&f_k-f^\ast-\sigma_1\alpha_1m^2\|x_k-x^\ast\|^2\cos^2\theta_k\\
&\leq&f_k-f^\ast-\alpha_2(f_k-f^\ast)\cos^2\theta_k\\
&=&(1-\alpha_2\cos^2\theta_k)(f_k-f^\ast)
\end{array}$$
for all sufficiently large $k$.\\

Hence, there is a constant $t=1-\alpha_2\cos^2\theta_k\in [0,1)$ such that
$$f(x_k)-f(x^\ast)\leq t^k[f(x_0)-f(x^\ast)].\eqno(5.9)$$\\

By Eq. 5.6, we have
$$f_k-f^\ast\geq\frac{1}{2}m\|x_k-x^\ast\|^2,$$
which, together with Eq. 5.9, implies
$$\|x_k-x^\ast\|\leq t^{k/2}[\frac{2(f(x_0)-f(x^\ast)}{m}]^{1/2}.$$
This is to say the sequence $\{x_k\}$ is $R$-linearly convergent. The proof of this theorem is complete.$\hfill \blacksquare$

\end{document}